\documentclass{amsart}

\usepackage{amsmath}
\usepackage{amsfonts}
\usepackage{amssymb}
\usepackage{mathrsfs} 
\usepackage{enumerate} 
\usepackage{tikz}
\usetikzlibrary{babel}
\usetikzlibrary{patterns}
\usepackage{tikz-cd}

\theoremstyle{plain}
\newtheorem{theo}{Theorem}[section]
\newtheorem{prop}[theo]{Proposition}

\newcommand{\To}{\Rightarrow}

\newcommand{\eps}{\varepsilon}
\newcommand{\N}{\mathbb N}
\newcommand{\Z}{\mathbb Z}
\newcommand{\R}{\mathbb R}

\newcommand{\one}{I}

\renewcommand{\AA}{\mathscr A}
\newcommand{\NN}{\mathscr N}
\newcommand{\BB}{\mathscr B}

\newcommand{\DD}{\mathscr D}
\newcommand{\EE}{\mathscr E}
\newcommand{\FF}{\mathscr F}
\newcommand{\UU}{\mathscr U}
\newcommand{\PP}{\mathscr P}

\renewcommand{\SS}{\mathscr S}

\newcommand{\leer}[1]{}

\newcommand{\ol}{\overline}

\newcommand{\sfrac}[2]{{\textstyle \frac{#1}{#2}}}
\newcommand{\symdif}{\vartriangle}

\DeclareMathOperator{\orb}{orb}
\DeclareMathOperator{\Int}{Int}
\DeclareMathOperator{\Aut}{Aut}
\DeclareMathOperator{\dist}{dist}

\newcommand{\ignore}[1]{}
\newcommand{\neu}[1]{{\sc #1}}
\newcommand{\salg}{$\sigma$-algebra }

\theoremstyle{plain}

\author{J. Wengenroth}
\address{Universit\"at Trier, FB IV -- Mathematik, 54286 Trier, Germany}
\email{wengenroth@uni-trier.de}
\title[Effros' theorem and descriptive set theory]{Effros' theorem on transitive group actions with a glimpse into descriptive set theory}

\begin{document}

 \begin{abstract}
The main aim of this note is to prove a version of a celebrated theorem of Effros \cite{Effros} about transitive group actions in a non-metrizable setting, these parts have been
formalized and verified with Lean by Lara Toledano.
We do not claim any originality since the given proof is in fact very close to one of van Mill \cite{vanMill}. Our presentation is however completely self-contained
and may serve as an appetizer to descriptive set theory. It also contains a few results about Suslin spaces (continuous images of
separable completely metrizable spaces, which are often called analytic) which are extremely useful in measure theory.
\end{abstract}

\maketitle

\section{Introduction}
The present notes are part of a bigger yet unfinished book project. The main reason for the separate publication is that the author's student Lara Toledano just finished a research
project in which she formalized the present proof of Effros' theorem and verified it with Lean. As far as we know this is the first formal verification of this important theorem.

Descriptive set theory is often presented in the realm of metrizable topological spaces which excludes many spaces appearing in analysis, e.g., Schwartz' distribution spaces
$\DD'(\Omega)$ of all distributions in an open set $\Omega\subseteq\R^d$, $\EE'(\Omega)$ of distributions with compact support, or $\SS'(\R^d)$ of tempered distributions, all of them endowed
with natural locally convex vector space topologies which are non-metrizable. 
It is thus desirable to avoid metrizability assumptions as far as possible -- and it turns out that very much is possible by only small adjustments of standard arguments.

\section{Suslin spaces}

A Hausdorff topological space $Y$ is called a \neu{Suslin space} (another transliteration of the Russian name is {\it Souslin}) if there are a polish space $X$ (i.e., a separable and complete metric space) 
and a continuous surjection $X\to Y$.

Trivially, every Suslin space is separable (even more, there is a countable {\it sequentially dense} subset)
and every coarser Hausdorff topology on a Suslin space makes the space
again Suslin.  A subset of a topological space which is Suslin in the relative topology is also  
called \neu{analytic} (sometimes, it is required that the superspace is metrizable), but we prefer the term Suslin subspace.
 There was an intense priority debate mixed with political accusations against Luzin about the contributions of Alexandrov, Luzin, Hausdorff and Suslin.
A detailed description  is in Lorentz' article \cite{Lorentz}.

\medskip To use a formulation of Schwartz \cite[page 116]{Schwartz},
there is a {\it tremendous stability of the Suslin spaces under the elementary operations on topological spaces} 
which yields that virtually all separable spaces appearing in analysis are Suslin. The Borel \salg of a topological space
$X=(X,\tau)$ is denoted by $\BB(X)$ or $\BB(\tau)$, it is the smallest \salg on $X$ which contains all open sets.
Elements of $\BB(X)$ are also called Borel subsets of $X$.

The topological \neu{coproduct} of a family of topological spaces $Y_\alpha$ is the disjoint (or {\it disjointified} -- although this word does not exist) union $Y=\coprod Y_\alpha$ endowed with the
topology consisting of all subsets $A$ such that all preimages $i_\alpha^{-1}(A)$ are open in $Y_\alpha$ where $i_\alpha:Y_\alpha\to Y$ denotes the inclusion of $Y_\alpha$ into the disjoint union.

\begin{theo}[Suslin stability]\label{thm:Suslinstability}
The class of Suslin spaces is stable under countable products, closed or open subspaces, countable topological coproducts, and continuous images with values in Hausdorff spaces. In particular,
countable intersections and unions of Suslin subspaces of a Hausdorff topological space are Suslin.
Every Borel subset of a Suslin space is Suslin. 
\end{theo}

\begin{proof}
Countable products or coproducts of polish spaces $(X_n,d_n)$ are again polish: The product $\prod X_n$ is endowed with the metric
\[
d((x_n)_{n\in\N},(y_n)_{n\in\N})=\sup \{d_n(x_n,y_n)\wedge 1/n: n\in\N\}
\]
(with $a\wedge b=\min\{a,b\}$) and the coproduct  $\coprod  X_n$  gets the distance 
\[
d(x,y)=d_n(x,y)\wedge 1\text{ if $x,y$ both belong to $X_n$ and $1$ otherwise.}
\]
For the coproduct, countability is only
needed to stay separable. If $Y_n=f_n(X_n)$ are continuous images, we get continuous surjections 
\[
\prod X_n\to\prod Y_n,\; (x_n)_{n\in\N}\mapsto (f_n(x_n))_{n\in\N}
\]
between the products and also a surjection $\coprod X_n \to\coprod Y_n$.

If $A$ is a closed or open subspace of a Suslin space $Y=f(X)$ for a continuous map $f:X\to Y$ on a polish space $(X,d)$, then $A=f(f^{-1}(A))$ is a continuous image
of $f^{-1}(A)$. If $A$ and hence $f^{-1}(A)$ are closed, this preimage is polish for the metric restricted from $X$, and if $A$ and hence $B=f^{-1}(A)$ are open, one can introduce in $B$
the complete metric
\[
d_B(x,y)=d(x,y) +|f(x)+f(y)| \text{ with } f(x)=1/\dist(x,X\setminus B)
\]
which penalizes closeness to the boundary of $B$ in order to become complete. In either case, $A$ is Suslin.

Countable intersections $\bigcap_{n\in\N} A_n$ of Suslin subspaces are homeomorphic to the (closed) diagonal in the product $\prod_{n\in\N} A_n$, and a countable union
$\bigcup_{n\in\N} A_n$ is a continuous image of the coproduct. Hence, the class of Suslin subspaces is stable with respect to countable intersections and unions.

To show that Borel subspaces of a Suslin space $Y$ are Suslin, we consider
\[
\AA=\{S\subseteq Y: S \text{ and } S^c \text{ are Suslin}\}.
\]
This class is trivially stable
under complements, and it is stable under countable unions (by de Morgan's rule). Therefore, $\AA$ is a \salg containing all open sets and hence all Borel sets.
\end{proof}

The following remarks explain our motivation to avoid metrizability in the treatment of Suslin spaces, they are not needed in the sequel.

In functional analysis (beyond Banach spaces) one often considers countable limits and colimits of locally convex spaces (whose topologies are generated by families of seminorms).
Limits in the locally convex category (with continuous linear maps as morphisms) are the same as those in the topological category, namely subspaces of products, and colimits (also called
inductive limits) are continuous images (even quotients) of coproducts which however are different from the coproducts in the topological category. The locally convex coproduct or, to
have different names in both categories, the \neu{direct sum} of locally convex spaces $X_\alpha$ is the vector space
\[
\bigoplus X_\alpha=\left\{(x_\alpha)_{\alpha\in I}\in\prod_{\alpha\in I} X_\alpha: \{\alpha\in I: x_\alpha\neq 0\} \text{ finite}\right\}
\]
endowed with the topology which is generated by all seminorms $p:\bigoplus X_\alpha\to [0,\infty)$ such that $p\circ j_\alpha$ are continuous
on $X_\alpha$ for all $\alpha$ where $j_\alpha: X_\alpha \to \bigoplus X_\alpha$ sends $x$ to the family $(x_\beta)_{\beta\in I}$ with $x_\alpha=x$ and $x_\beta=0$ for $\beta\neq \alpha$.
This topology is (usually much) coarser than the finest topology making all $j_\alpha$ continuous but it is Hausdorff if so are all $X_\alpha$. For a countable index set $I=\N$, we get
a continuous surjection $\coprod_n\left( \prod_{k\le n} X_k\right) \to \bigoplus_\N X_n$ so that countable direct sums of Suslin locally convex spaces are Suslin.

Since locally convex colimits are quotients and hence continuous images of direct sums, the class of Suslin spaces is stable with respect to countable
locally convex colimits.

\medskip
This discussion yields, e.g., that the space $\DD(\Omega)$ of {\it test functions} is a Suslin space as a Hausdorff countable colimit of the separable Fr\'echet spaces
$\DD(K_n)$ of smooth functions with support in a fixed compact subset $K_n$ of $\Omega$ with the topology of uniform topology of
all partial derivatives. It follows, e.g., from theorem \ref{thm:factorization groups} below that $\DD(\Omega)$ is not metrizable. More generally,
if a countable locally convex colimit $X$ of Fr\'echet spaces $X_n$ is metrizable and complete then $X=X_n$ for some $n\in\N$. 

The space $\DD'(\Omega)$ of distributions in $\Omega$ is a countable limit of a countable colimit of separable Banach spaces
and hence also Suslin. The book \cite[page 115]{Schwartz} of Schwartz contains many more example of spaces related to partial differential equations and Fourier transformation which are Suslin.
As a rule of thumb, all separable spaces one usually meets in analysis are Suslin but many of them are not metrizable.

\section{The Borel \salg of Suslin spaces}

A very remarkable property of Suslin spaces is the behaviour of the Borel \salg in the following theorem:

\begin{theo}[Luzin, Suslin]\label{thm:Luzin}
If $(X,\tau)$ is a Suslin space and $\tilde\tau$ is a coarser Hausdorff topology on $X$ then $\BB(\tau)=\BB(\tilde\tau)$.
\end{theo}

Here is a typical application:

\medskip
{\it Let $\varphi_n:X\to Y_n$ be a point separating sequence of continuous functions on a Suslin space $X$ with values in Hausdorff spaces $Y_n$ which have
countable bases $\UU_n$ of their topologies (e.g., separable metric spaces). For a measurable space $(\Omega,\AA)$, a function
$f:\Omega\to X$ is Borel measurable if (and only if) all $\varphi_n\circ f:\Omega\to Y_n$ are Borel measurable. 
}

\medskip
This follows from the theorem for the initial topology $\tilde\tau$ of $\{\varphi_n:n\in\N\}$ which is coarser than $\tau$ because all $\varphi_n$ are continuous
and which is Hausdorff because the sequence separates the points of $X$. Moreover, $\tilde\tau$
has a countable basis $\UU$ consisting of all finite intersections of preimages $\varphi_n^{-1}(B_n)$ with $B_n\in\UU_n$. Since every $\tilde \tau$-open set is a {\it countable} union of elements from $\UU$, we get
$\BB(\tau)=\BB(\tilde\tau)=\sigma(\UU)$ and the measurability of all $\varphi_n\circ f$ implies that $f$ is $(\AA,\sigma(\UU))$-measurable.

\medskip

The proof of the theorem needs some preparation. % (our  presentation  follows Schwartz' book {\it Radon Measures on Arbitrary Topological Spaces and Cylindrical Measure}). 
We will first show that every Suslin space is a continuous image of $\N^\N$, i.e., the countable product of the integers (where $\N$ is endowed with the discrete topology).

For $\alpha=(\alpha_k)_{k\in\N} \in\N^\N$ and $n\in\N_0$ we denote by $\alpha|n=(\alpha_1,\ldots,\alpha_n)$ the {\it truncation at $n$} (in particular, $\alpha|0$ is the
empty string $(\,)$).
A \neu{Suslin representation} of a Hausdorff topological space $X$ is a family
$\FF=\{F_{\alpha|n}: \alpha\in\N^\N, n\in\N\}$ of subsets of $X$ with the following properties:
\begin{enumerate}
\item[(I)] $X=F_{\alpha|0}$ and $F_{\alpha|n}=\bigcup \{F_{\beta|n+1}: \beta\in\N^\N, \beta|n=\alpha|n\}$
for all $\alpha\in\N^\N$ and $n\in\N_0$.
\item[(II)]
For every $\alpha\in \N^\N$ the sequence  % bei Schwartz sogar ohne Abschluss \bigcap_{N\in\N} \overline{F_{\alpha|N}}$
 $(F_{\alpha|n})_{n\in\N}$ converges to some $x\in \bigcap_{n\in\N} F_{\alpha|n}$, i.e., every neighbourhood of $x$
contains $F_{\alpha|n}$ for $n\in\N$ big enough.
\end{enumerate} 

It should be noted that despite the uncountability of $\N^\N$, a Suslin representation $\FF$ is a countable system because it can be indexed by $\bigcup_{n\in\N_0} \N^n$.
The first condition can also be written as $X=\bigcup_{k\in\N} F_{(k)}$ and $F_{(\alpha_1,\ldots,\alpha_n)}=\bigcup_{k\in\N}  F_{(\alpha_1,\ldots,\alpha_n,k)}$. Only the second condition depends
on the topology of $X$, Hausdorffness implies that the limit $x\in X$ of $F_{\alpha|n}$ is unique, the requirement $x\in \bigcap_{n\in\N} F_{\alpha|n}$ is automatically fulfilled
if all $F_{\alpha|n}$ are closed. 

If $\FF$ is a Suslin representation of $X$ and $f:X\to Y$ is a continuous map onto a Hausdorff space $Y$ then
$\{f(F_{\alpha|n}): \alpha\in\N^\N, n\in\N_0\}$ is a Suslin representation of $Y$.

\begin{theo}[Suslin representations]\label{thm:Suslin representation}

For a non-empty Hausdorff topological space $X$, the following conditions are equivalent.
\begin{enumerate}
\item[(1)] $X$ is Suslin,
\item[(2)] $X$ has a Suslin representation,
\item[(3)] $X$ is a continuous image of $\N^\N$.
\end{enumerate}
\end{theo}

The implication (2) $\To$ (1) shows that the elements of a Suslin representation $\FF$ are themselves Suslin spaces. 

\begin{proof}
(2) $\To$ (3).
If $\FF$ is a Suslin representation of $X$, we get a map $f:\N^\N\to X$ assigning to $\alpha\in\N^\N$  the limit of the sets $F_{\alpha|n}$. Since $X$ is Hausdorff, this limit is unique
so that $f$ is indeed a function into $X$. Let us show that $f$ is continuous. For $\alpha\in\N^\N$
and a neighbourhood $U$ of $x=f(\alpha)$ there is $n\in\N$ with $F_{\alpha|n}\subseteq U$. Then $V=\{\beta\in\N^\N: \beta|n=\alpha|n\}$ is a neighbourhood of $\alpha$ with
$f(\beta)\in \bigcap_{k\in\N} F_{\beta|k} \subseteq F_{\beta|n}=F_{\alpha|n}\subseteq U$ for all $\beta\in V$.

Surjectivity of $f$ follows from the first condition of Suslin representations: For $x\in X$, there is $\alpha_1\in\N$ with $x\in F_{(\alpha_1)}=\bigcup_{\alpha_2\in\N} F_{(\alpha_1,\alpha_2)}$, hence there is
$\alpha_2\in\N$ with $x\in F_{(\alpha_1,\alpha_2)}=\bigcup_{\alpha_3\in\N} F_{(\alpha_1,\alpha_2,\alpha_3)}$. Recursively, we thus find $\alpha\in\N^\N$
with $x\in\bigcap_{n\in\N} F_{\alpha|n}$ and hence $f(\alpha)=x$.

\medskip
(3) $\To$ (1) is clear since $\N^\N$ is polish but it may be worthwhile to note that a particularly simple Suslin representation of $\N^\N$ is given by
$F_{\alpha|n}=\{\beta\in\N^\N: \beta|n=\alpha|n\}$.

For the proof of (1) $\To$ (2), it is enough to construct a Suslin representation for every polish space $(X,d)$.
For $\eps=1/2$ and a countable dense set $S$ of $X$ the family of open balls of radius $\eps$ and centre in $S$ is an open cover $\{E_{\alpha_1}:\alpha_1\in\N\}$ of $X$ 
(we do not require the sets to be distinct, it may happen that we have to repeat a single set infinitely many times). We set $F_{(\alpha_1)}=\overline{E}_{\alpha_1}$
(which has diameter $\le 1$) 
and apply the same argument  to the cover of $F_{(\alpha_1)}$ by countably many open balls of radius $1/4$. This
yields an open cover $\{E_{\alpha_1,\alpha_2}:\alpha_2\in\N\}$ of $F_{(\alpha_1)}$ and we set $F_{(\alpha_1,\alpha_2)}=F_{(\alpha_1)}\cap
\overline{E}_{\alpha_1,\alpha_2}$ which has diameter $\le 1/2$. Continuing in this way we thus get a
family $\{F_{\alpha|n}: \alpha\in\N^\N, n\in\N\}$ of non-empty closed sets $F_{\alpha|n}$ of diameter $\le 1/2^{n-1}$.  Condition (I) is fulfilled by construction, and condition
(II) is a consequence of the completeness: The intersection of every decreasing sequence of non-empty closed sets with diameters tending to $0$ is a singleton and a limit
of the sequence.
\end{proof}

\pagebreak[3]

\begin{theo}[Luzin's separation theorem]\label{thm:Luzin separation}
 Any  two disjoint Suslin subsets $A$ and $\tilde A$ of a Hausdorff topological space $X$ are \neu{Borel separated}, i.e., there are disjoint Borel sets
$B$ and $\tilde B$ of $X$ with $A\subseteq B$ and $\tilde A \subseteq \tilde B$ (equivalently, $A\subseteq C \subseteq  X\setminus \tilde A$ for a Borel set $C$). In particular, if $A$ and $A^c$ are both Suslin then
$A$ is Borel.
\end{theo}

\begin{proof} We will use the following principle recursively:

\medskip
{\it
If $E_n$ and $\tilde E_n$ are sets with $A=\bigcup_{n\in\N} E_n$ and $\tilde A=\bigcup_{n\in\N} \tilde E_n$ such that, for all $n,m\in\N$, $E_n$ and $\tilde E_m$ are separated by a Borel set $B_{n,m}$,
i.e.,  $E_n\subseteq B_{n,m}\subseteq \tilde E_m^c$,
then the Borel set $B=\bigcup_{n\in\N}\bigcap_{m\in\N} B_{n,m}$ separates $A$ and $\tilde A$.}

\medskip
Let $\FF$ and $\tilde \FF$ be Suslin representations of $A$ and $\tilde A$, respectively. Assuming that $A$ and $\tilde A$ are not Borel separated we apply the principle above
to $A=\bigcup_{n\in\N} F_{(n)}$ and $\tilde A=\bigcup_{n\in\N} \tilde F_{(n)}$. We thus find $\alpha_1\in\N$ and $\tilde\alpha_1\in\N$ such that 
$F_{(\alpha_1)}$ and $\tilde F_{(\tilde \alpha_1)}$ are not Borel separated. Using the principle again for $F_{(\alpha_1)}=\bigcup_{n\in\N} F_{(\alpha_1,n)}$ and the corresponding
cover of $F_{(\tilde\alpha_1)}$, we find $\alpha_2\in\N$ and $\tilde\alpha_2\in\N$ such that $F_{(\alpha_1,\alpha_2)}$ and $\tilde F_{(\tilde\alpha_1,\tilde\alpha_2)}$ are not Borel separated.
Continuing in this way, we find $\alpha,\tilde\alpha \in\N^\N$ such that $F_{\alpha|n}$ and $\tilde F_{\tilde\alpha|n}$ are not Borel separated for every $n\in\N$.
These sets converge to $x\in A$ and $\tilde x\in\tilde A$, respectively. Since $A$ and $\tilde A$ are disjoint and $X$ is Hausdorff there are disjoint open  neighbourhoods $U$ of $x$ and 
$\tilde U$ of $\tilde x$. The convergence of the sets thus yields $n\in\N$ with $F_{\alpha|n}\subseteq U$ and $\tilde F_{\tilde\alpha|n}\subseteq \tilde U$. Since open sets
are Borel we have thus found a Borel separation, contradicting the construction of $\alpha$ and $\tilde\alpha$.
\end{proof}

\begin{proof}[Proof of theorem \ref{thm:Luzin}] From $\tilde\tau\subseteq\tau$ we get $\BB(\tilde\tau)\subseteq \BB(\tau)$. If, on the other hand,
$A$ is a Borel subset of $(X,\tau)$ then $A$ and $A^c$ are both $\tau$-Suslin since Borel subsets of Suslin spaces are Suslin by theorem \ref{thm:Suslinstability}. 
Continuity of the identity map $(X,\tau)\to (X,\tilde\tau)$ implies that $A$ and $A^c$ are both $\tilde\tau$-Suslin and hence
$\tilde\tau$-Borel by Luzin's separation theorem.
\end{proof}

For an application of the separation theorem, we need another property (which plays a much bigger role in classical descriptive set theory than in our
little introduction). A subset $B$ of a topological space has the \neu{Baire property} if it differs from some open set
$U$ only by a meagre set $M$, i.e., $B \symdif U =M$ with the symmetric difference $B\symdif U= (B\setminus U) \cup (U\setminus B)$.
Recall that a subset $M$ of a topological space is \neu{meagre} (the reader may have noticed the author's liking for category theory so that
we avoid the alternative terminology {\it first category sets}) if $M$ is contained in a countable union of closed sets with empty interiors.

\medskip

The power set of $X$ is a commutative ring for the operations $\symdif$ as the sum (with the peculiar property
$A\symdif A=\emptyset$ for all $A$) and $\cap$ as the multiplication.
This is convenient for calculations, for example, $B\symdif U=M$ implies 
\[
B=B\symdif (U \symdif U)=(B\symdif U) \symdif U=M\symdif U
\]
so that sets with Baire property are symmetric differences of an open and a meagre set.
The boundary $\partial U = \overline U \setminus U$ of an open set $U$ is closed with empty interior and hence meagre. Since $U$ and $\partial U$ are disjoint
we then obtain
\[
B\symdif \overline U = B\symdif (U\symdif \partial U)=(B\symdif U)\symdif \partial U \subseteq (B\symdif U) \cup \partial U.
\]
 This yields that every set with the Baire property also differs from some closed $F$ set by a meagre set. Conversely, the latter condition implies the
 Baire property  since, for the interior $U$ of $F$, we have $B\symdif F =(B\symdif U) \symdif \partial U$ which implies
 $B\symdif U = (B\symdif F)\symdif \partial U \subseteq (B\symdif F) \cup \partial U$. 
 
 Since $B^c\symdif U^c=B\symdif U$ we obtain
 that the system $\BB\PP(X)$ of subsets with the Baire property is stable under complements. From
 \[
 \left( \bigcup_{n\in\N} B_n\right) \symdif \left(\bigcup_{n\in\N} U_n\right) \subseteq \bigcup_{n\in\N} B_n \symdif U_n
 \]
 we also get that $\BB\PP(X)$ is stable under countable unions so that $\BB\PP(X)$ is a \salg which contains all open sets (since $U\symdif U=\emptyset$) as well as all meagre sets
 (since $M\symdif\emptyset=M$).  We would like to warn the reader that the name {\it Baire \salg} is used in for the smallest \salg making all real valued continuous functions measurable,
 this \salg is smaller than the Borel \salg whereas $\BB\PP(X)$ is typically strictly bigger than the Borel $\sigma$-algebra.
 
 \medskip
 We mention a rather deep result of Suslin that Suslin subsets of separable metrizable spaces also have the Baire property, but we will only
 use the much simpler fact just proved that Borel sets in Hausdorff spaces have the Baire property. It turns out to be quite
 advantageous that we do not need metrizability in the following theorem of van Mill which however in \cite{vanMill} is stated only for
 metrizable spaces.

\begin{theo}[van Mill]\label{thm:van Mill}
Let $S,T$ be two Suslin subsets of a Hausdorff topological space such that $S$ is not meagre and $T$ is \neu{nowhere meagre}, i.e., the intersection with every non-empty open set is
non-meagre. Then $S\cap T$ is not empty.
\end{theo}
 
 Without the Suslin property of $S$ and $T$, the result is wrong even for subsets of $\R$, examples are {\it Bernstein sets} $S\subseteq \R$ such that $S$ and $T=S^c$ meet
 every closed uncountable subset of $\R$. A construction of such sets (using a well-order of $\R$) can be seen
 in theorem 5.3 of Oxtoby's book \cite{Oxtoby} and lemma 5.1 there shows that $S$ and $T$ are nowhere meagre.

 \begin{proof}
Assuming $S\cap T=\emptyset$ the separation theorem yields a Borel set $B$ with $S\subseteq B\subseteq X\setminus T$. In particular, $B$ is again non-meagre
and we will reach a contradiction by showing $B\cap T\neq \emptyset$. Since $B$ has the Baire property there are an open set $U$ and a meagre set $M$ with $B=U\symdif M$.
The set $U$ is not empty because $B$ is not meagre. Hence $B\cap T= (U\cap T) \symdif (M\cap T) \supseteq (U\cap T) \setminus M$ is not meagre since so is $U\cap T$.
\end{proof}

 The result just proved is one of the main ingredients for the open mapping principle in Effros' theorem \ref{thm:Effros} and its consequences. 
 Before turning to this we mention that the version of theorem \ref{thm:van Mill} given here is not optimal.
 In separable metrizable spaces $X$, one can use the Baire property of Suslin subsets to show that the
 intersection $S\cap T$ is even non-meagre. Moreover, one can extend the theorem to countable intersections of nowhere meagre Suslin sets in separable metrizable spaces.
Readers who like to dive much deeper into descriptive set theory are referred to the book \cite{Kechris} of Kechris.

\section{Effros' theorem}
\medskip
We now come to \neu{group actions} $G\times X\to X$, $(g,x)\mapsto g.x$ of a topological group $G$ on a topological space $X$. Algebraically, this
means $g.(h.x)=(gh).x$ for all $g,h\in G$ (the group operation is written by juxtaposition) and $e.x=x$ for the neural element $e$ of $G$ and  all $x\in X$.
That $G$ is a topological group means that $G$ is endowed with a topology such that $G\times G\to G$, $(x,y)\mapsto xy^{-1}$ is continuous.

 For
$x\in X$, the set $\orb(x)=\{g.x: g\in G\}$ is the \neu{orbit} of $x$ under the action. The action is called \neu{transitive} if there is
only a single orbit or, equivalently, for all $x,y\in X$, there is $g\in G$ with $g.x=y$. One then wants to know whether, for $y\in X$ sufficiently close to
$x\in X$ , one can {\it move} $x$ to $y=g.x$ with $g$ close to the identity, such group actions are called \neu{micro-transitive}.

The case of a polish group acting transitively on a polish space in the next result is the classical theorem of Effros \cite{Effros}, the
version given here is essentially the one of van Mill \cite{vanMill} who considered transitive actions of metrizable Suslin groups.

\begin{theo}[Effros]\label{thm:Effros}
Let $G$ be a topological group whose topology is Suslin, $X$ a Hausdorff topological space, and $G\times X\to X$, $(g,x)\mapsto g.x$ a separately continuous group
action such that no orbit is meagre in $X$. Then $G\to X$, $g\mapsto g.x$ is open for every $x\in X$ and all orbits are closed and open.
\end{theo}

 \begin{proof}
We will first show that $U.x$ is not meagre for every $x\in X$ and $U\in\NN_e(G)$. For a countable dense set $\{g_n:n\in\N\}$ of $G$, we have $G=\bigcup_{n\in\N} g_nU$ (for $g\in G$, the neighbourhood
$gU^{-1}$ of $g$ intersects the dense set) so that $\orb(x)=\bigcup_{n\in\N} g_nU.x= \bigcup_{n\in\N} g_n.(U.x)$. Since $\orb(x)$ is not meagre there is $n\in\N$ such that $g_n.(U.x)$ is not meagre and since
translation by $g_n$ is a homeomorphism also $U.x$ is not meagre.

\medskip
The next step is to show that, for all $x\in X$, open $U\in\NN_e(G)$, and open sets $A\subseteq X$, the intersection $A\cap U.x$ is either empty or non-meagre.

Indeed, for $a\in A\cap U.x$ and $u\in U$ with $a=u.x$, the set $V=\{g\in G: g.x\in A\}$ is an open neighbourhood of $u$ (by the continuity of $g\mapsto g.x$) so that
$W=Vu^{-1} \cap Uu^{-1} \in\NN_e(G)$. For $w\in W$, we then have $w.a=wu.x\in V.x\cap U.x\subseteq A\cap U.x$. Since $W.a$ is not meagre by the first step so is $A\cap U.x$.

\medskip
For $x\in X$ and open $U\in \NN_e(G)$, we will next show that the interior $B=\Int\left(\ol{U.x}\right)$ is dense in $\ol{U.x}$ and satisfies $x\in B$.

Indeed, we have to show that $B$ meets every open neighbourhood $A$ of any given element of $\ol{U.x}$. Since $M=A\cap U.x\neq \emptyset$, this set is not meagre
so that $\Int(\ol M)$ has an element $y$. Then $B\supseteq \Int(\ol M)$ is a neighbourhood of $y\in \ol A$ so that $B\cap A\neq\emptyset$.
To show $x\in B$ we choose $V\in\NN_e(G)$ with $V^{-1}V\subseteq U$. The first step implies that $C=\Int(\ol{V.x})$ is not empty and since $C$ is a neighbourhood of each of its elements
we have $C\cap V.x\neq\emptyset$. There is thus $v\in V$ with $v.x\in C$ and using that the translation by $v^{-1}$ is a homeomorphism we obtain
\begin{align*}
x & =v^{-1}.(v.x)\in v^{-1}.C =v^{-1}.\Int\left(\ol{V.x}\right)=\Int\left(\ol{v^{-1}.(V.x)}\right) \\ & \subseteq \Int\left(\ol{V^{-1}V.x}\right) \subseteq \Int\left(\ol{U.x}\right)=B.
\end{align*}

\medskip
For open $U,V\in\NN_e(G)$ with $V^{-1}V\subseteq U$ and $x\in X$, we will finally prove that $A=\Int(\ol{V.x})\subseteq U.x$ which by the previous step shows that
$U.x$ is a neighbourhood of $x$ so that $G\to X$, $g\mapsto g.x$ is open at $e$. Translation invariance then yields that the $g\mapsto g.x$ is open.

Given $y\in A$ we will prove $V.x \cap V.y\neq \emptyset$. If then $v.x=w.y$ with $v,w\in V$ we get $y=w^{-1}v.x\in U.x$, as required.
We set $B=\Int(\ol{V.y})$ so that, by the previous step, $C=A\cap B$ is an open neighbourhood of $y$.   The set $V.x$ is the continuous image of the
open subset $V$ of the Suslin space $X$ so that $V.x$ is again Suslin and hence so is the relatively open subset $V.x \cap C$
by theorem \ref{thm:Suslinstability}. In order to apply van Mill's theorem \ref{thm:van Mill} we show that $V.x\cap C$ is {\it nowhere} meagre in $C$: If $D$ is an open ($=$ relatively open)
non-empty subset of $C$ we get $V.x\cap D\neq \emptyset$ because $C$ is contained in $\ol{V.x}$. The second step of the proof thus yields
that $V.x\cap D$ is non-meagre in $X$ and thus also non-meagre in $C$. The same argument shows that
$V.y\cap C$ is a nowhere meagre Suslin subset of $C$ so that theorem \ref{thm:van Mill} yields $V.x\cap V.y\neq \emptyset$.

\medskip

It remains to show that $G.x$ is closed for every $x\in X$. For $y\in\ol{G.x}$, the orbit $G.y$ is open and hence a neighbourhood of $y$ which thus intersects $G.x$. For $g,h\in G$ with
$g.y=h.x$ we obtain $y=g^{-1}h.x\in G.x$.
\end{proof}

Let us briefly consider some  topological-geometric background of Effros' theorem.
A topological space $X$ is called \neu{homogeneous} if the group $\Aut(X)$ acts transitively on $X$, i.e., for all $x,y\in X$, there is an automorphism $g$ with $g(x)=y$. Since
neighbourhoods of $x$ and $y$ are then homeomorphic, homogeneity implies that {\it the space looks the same everywhere} (this latter property however is local and does not
imply homogeneity, the disjoint union of a circle and an open interval is an example).

The class of homogeneous spaces
is stable with respect to homeomorphisms and with respect to products but has otherwise rather poor stability properties, the interval $[0,1]$
is not homogenous because small neighbourhoods of the boundary and interior points are not homeomorphic. Topological groups are obviously
homogeneous since the translations are automorphisms (so that also the Cantor set as well as $\N^\N$ are homogeneous because they are
homeomorphic to $\{0,1\}^\N$ and $\Z^\N$, respectively) 
and  so are many spaces which are directly related to some group action, e.g., the left or right quotients of a topological group with respect
not necessarily normal subgroup. As non-obvious examples  we mention that the {\it Hilbert cube} $[0,1]^\N$ is
homogeneous (this is proved as theorem 6.1.6 in van Mill's book \cite{vanMillBook})   %{\it Infinite-Dimensional Topology: Prerequisites and Introduction}) 
and the following proposition (which is certainly known although it is hard to find a precise reference).

\begin{prop}
Every connected Hausdorff space $X$ which is {\it locally homeomorphic to open subsets of normed spaces}, e.g., a connected
topological manifold, is homogeneous.
\end{prop}

\begin{proof}
Every point $x\in X$ has a neighbourhood basis of sets homeomorphic
to the closed unit ball $B$ of a normed space. For any point $p\in \Int(B)$, there is an automorphism of $B$ moving the origin to $p$ and fixing the
boundary, e.g., $\varphi(x)=x+(1-\|x\|)p$ which maps each {\it radial segment} from the origin to $y\in\partial B$ to the segment $[p,y]$. For the inverse, 
take $y\in B\setminus\{p\}$ and the {\it unique} $s(y)\ge 1$ such that $z(y)=p+s(y)(y-p)$ has norm $1$  (uniqueness follows from the convexity
of $f(s)=\|p+s(y-p)\|$ which has an increasing right derivative and satisfies $f(0)<1$), then continuity of $s(y)$ and hence $\varphi^{-1}(y)=\sfrac{1}{s(y)}z(y)$ can be seen
with a sub-subsequence argument. 

Using this for
slightly smaller closed neighbourhoods and extending the constructed automorphism by the identity on the complement of that neighbourhood,
one gets a neighbourhood basis of $x$ consisting of open sets $U$ such that, for every $y\in U$ there is  an automorphism of $X$ moving
$x$ to $y$. For fixed $x_0\in X$, this implies that $A=\{x\in X: f(x_0)=x \text{ for some $f\in \Aut(X)$}\}$ is open and
also that $X\setminus A$ is open. Connectedness thus implies $A=X$.
\end{proof}

\medskip

For a compact metric space $(X,d)$, the automorphism group has a natural metric $d_\infty(f,g)=\sup\{d(f(x),g(x)):x\in X\}$ which is clearly right invariant and makes $\Aut(X)$ a topological group
(the elementary proof uses the {\it uniform} continuity of continuous functions on compact metric spaces).  Continuity of the inversion yields that
$d_\infty(f,g)\vee d_\infty(f^{-1},g^{-1})$ generates the same topology and since this metric is complete, $\Aut(X)$ is completely metrizable. 
To prove separability of $\Aut(X)$ it is enough to show that the  space $C(X,X)$ of all continuous functions with the metric $d_\infty$ of uniform convergence is separable.
Separability of the Banach space $C(X,\R)$ follows, e.g., from the Stone-Weierstra{\ss} theorem and since $X$ is homeomorphic to a compact subset
$[0,1]^\N$ (for a countable dense set $\{s_n:n\in\N\}$ of $X$, the map $X\to [0,1]^\N$, $x\mapsto (d(x,s_n)\wedge 1)_{n\in\N}$ is a continuous injection and hence a homeomorphism onto its image)
yields that $C(X,X)$ is homeomorphic to a closed subset of $C(X,\R)^\N$ and hence separable.

Therefore, $\Aut(X)$ is a polish group (this is useful also in other interesting cases, e.g., the group $G$ of isometric isomorphisms of a Banach
space with separable dual endowed with the operator norm is polish because it embeds into $\Aut(X)$ by assigning to $T\in G$
the restriction of the bitransposed to the unit ball of the bidual with the weak$^*$-topology which is compact by Alaoglu's theorem 
and metrizable since the Banach space has a separable dual).  Since $\Aut(X)$ acts continuously on $X$ by evaluation Effros' theorem implies:

\begin{theo}[Micro-transitivity]\label{thm:microtransitive}
For every homogeneous compact metric space $(X,d)$ and $\eps>0$, there is $\delta>0$ such that, for all $x,y\in X$ with $d(x,y)<\delta$, there is an automorphism
$f\in\Aut(X)$ with $f(x)=y$ and $\sup\{d(f(t),t): t\in X\} <\eps$.
\end{theo}
 
 \medskip

\section{Open mapping theorems}

Effros' theorem immediately yields a very general open mapping theorem for topological groups because
every group morphism $f:X\to Y$ between groups induces a group action $X\times Y\to Y$, $(x,y)\mapsto f(x)y$ whose
orbits are $\orb(y)=f(X)y$. If $Y$ is a topological group all orbits are
homeomorphic to the image $f(X)$ (this only needs separate continuity of the group operation and existence of inverses but
not the continuity of the inversion) so that we obtain an immediate corollary of theorem \ref{thm:Effros}.

\begin{theo}(Suslin open mapping theorem)\label{thm:SuslinOMT}
Let $X$ be a Suslin topological group, $Y$ a Hausdorff topological group, and $f:X\to Y$ a continuous group morphism such that
$f(X)$ is not meagre in $Y$. Then $f$ is open and $f(X)$ is closed and open in $Y$.
\end{theo}

The hypothesis that $f$ has non-meagre image is satisfied if $f$ is surjective and $Y$ is a Baire space. On the other hand, surjectivity of $f$ is a consequence
of the theorem if $Y$ is connected. But even if one is only interested in connected topological groups or, as a particular case, in topological vector spaces it can be
advantageous to have the present version of the theorem. An example is the following theorem 
similar to Grothendieck's factorization theorem \cite[Th\'eor\`eme A]{Grothendieck} for Fr\'echet spaces  (which does not need separability).

\begin{theo}[Grothendieck's factorization theorem for groups]\label{thm:factorization groups}
Let $f:X \to Y$ be a continuous group morphism from a polish topological group to a Hausdorff topological group and $g_n:Y_n\to Y$ continuous
group morphisms on Suslin topological groups  $Y_n$ such that $f(X)\subseteq \bigcup_{n\in\N} g_n(Y_n)$. Then there is $n\in\N$ such that, 
for every $U\in\NN_e(Y_n)$, there is $V\in \NN_e(X)$ with $f(V)\subseteq g_n(U)$. If $X$ is connected also $f(X)\subseteq g_n(Y_n)$ holds and if,
additionally, $g_n$ are injective there
is a continuous group morphism $h:X\to Y_n$ with $f=g_n\circ h$ as in the following diagram
\[
\begin{tikzcd}
X \arrow[rr,"f"] \arrow[drr, dashed, bend right=70, "h"']  & & Y& \\
 & Y_1\arrow [ur,near start,"g_1" ] & Y_n\arrow[u, near start,"g_n"] & Y_m. \arrow[ul,near start,"g_m"']
\end{tikzcd}
\]
\end{theo}

\begin{proof}
For  $n\in \N$, the {\it pull back} $Z_n=\{(x,y)\in X\times Y_n: f(x)=g_n(y)\}$ is closed in $X\times Y_n$ because of the continuity of $f$ and $g_n$ into the Hausdorff space $Y$. Hence
$Z_n$ is a Suslin topological group as a closed subgroup of $X\times Y_n$. The projections $\pi_n:Z_n\to X$, $(x,y)\mapsto x$ are continuous and the assumption $f(X)\subseteq \bigcup_{n\in\N} g_n(Y_n)$
implies $X=\bigcup_{n\in\N} \pi_n(Z_n)$. Baire's theorem implies that $\pi_n$ has non-meagre image for some $n\in\N$ and the open mapping theorem then
yields that $\pi_n$ is open. For every $U\in\NN_e(Y)$, there is thus $V\in\NN_e(X)$ with $V\subseteq \pi_n(X\times U)$ which implies $f(V)\subseteq g_n(U)$. 
Since $\pi_n(Z_n)$ is closed and open connectedness of $X$ implies $\pi_n(Z_n)=X$ and hence $f(X)\subseteq g_n(Y_n)$.
If $g_n$ is injective with inverse $g_n^{-1}: g_n(Y_n)\to Y_n$, the openness condition $f(V)\subseteq g_n(U)$ implies the continuity of $h=f\circ g_n^{-1}:X\to Y_n$ which is the desired
factorization $f=g_n\circ h$.  
\end{proof}

The open mapping theorem implies the following results which conclude continuity from
properties of the graph. 

\begin{theo}[Suslin and Borel graph theorem]\label{thm:Suslin-Borel-graph}
Let $X$ and $Y$ be Hausdorff topological groups  and $f:X\to Y$ be a group morphism.
\begin{enumerate}
\item If the graph $G(f)$ is a Suslin subspace of $X\times Y$ and $X$ is not meagre in itself then $f$ is continuous.
\item If  $X$ is completely metrizable, $Y$ is Suslin, and the graph of $f$ is Borel in $X\times Y$ then $f$ is continuous.
\end{enumerate}
\end{theo}

\begin{proof} The first part is proved as an application of the
open mapping theorem to the continuous bijective group morphism $g:G(f)\to X$, $(x,f(x))\mapsto x$ in the diagram
\[
\begin{tikzcd}
G(f)\arrow[r, "i"] \arrow[d, shift right=0.5ex, "g"']  \arrow[d, leftarrow, shift left=1ex, "g^{-1}"] & X\times Y  \arrow[d,"\pi_2"]  \\
X \arrow[r, "f"]                    & Y.
\end{tikzcd}
 \]
Then $g$ is open so that $g^{-1}$ and hence $f=\pi_2\circ i \circ g^{-1}$ are continuous.

\medskip
(b) It is enough to show sequential continuity and hence continuity of the restriction 
$\tilde f=f|_{\tilde X}:\tilde X\to Y$ for every separable closed subgroup of $X$.
 The graph $G(\tilde f)$ is the continuous preimage
of $G(f)$ under the continuous (and thus Borel measurable) inclusion $\tilde X\times Y\to X\times Y$ so that $G(\tilde f)$
is a Borel subset of the Suslin space $\tilde X\times Y$ and by theorem \ref{thm:Suslinstability} therefore itself Suslin. Part (a) then implies that $\tilde f$ is continuous.
\end{proof}

In linear functional analysis, de Wilde \cite{deWilde} proved general open mapping and closed graph theorems for the class of {\it webbed spaces} whose definition is somewhat technical.
His arguments are in the same spirit as the proof of Effros' theorem, the main difference is that the open mapping theorem for webbed spaces does not need separability.
The reason is that, for neighbourhoods $U$ of the origin of a topological vector space $X$, one always has a countable cover $X=\bigcup_{n\in\N} nU$ which uses the scaling
$\mathbb K\times X\to X$. Since this does not exist for general groups one has to use instead the cover $X=\bigcup_{n\in\N} g_nU$ for a dense set $\{g_n:n\in\N\}$.

\section{Measures on Suslin spaces}

We finish our glimpse into descriptive set theory with an important measure theoretic aspect of Suslin spaces. A measure
$\mu$ on the Borel \salg of a topological space is called \neu{tight} (a more expressive name would be {\it inner compactly
approximable}) if, for all Borel sets $A$,
\[
\mu(A)=\sup\{\mu(K): K\subseteq A \text{ compact}\}.
\]

The proof of the following theorem is a blend of topological and measure theoretic  ingredients and needs two little preparations.
For a measure space $(\Omega,\AA,\mu)$ and arbitrary subset $B$ of $\Omega$, we set
\[
\mu^*(B)=\inf\{\mu(A): B\subseteq A\in\AA\}
\]
which is the {\it outer measure} corresponding to $\mu$. Monotonicity of $\mu$ implies that
$\mu^*$ is again increasing and that $\mu^*(A)=\mu(A)$ for all $A\in \AA$. Moreover, we will
need that for every sequence $B_n\uparrow B$ we have $\mu^*(B)\le \lim\limits_{n\to\infty}\mu^*(B_n)$
(where this monotone limit exists in $[0,\infty]$). This is trivial if the right hand side is $\infty$
and otherwise we show the inequality up to an additive $\eps>0$. For each $n\in\N$, the definition of the
greatest lower bound yields
$A_n\in\AA$ with $B_n\subseteq A_n$ and $\mu(A_n)< \mu^*(B_n)+\eps$. The {\it limit inferior}
$A=\bigcup_{n\in\N}\bigcap_{k\ge n}A_k\in\AA$ then satisfies $B\subseteq A$ and Fatou's lemma implies
\begin{align*}
\mu^*(B) & \le \mu(A)=\int \one_Ad\mu =\int \liminf_{n\to\infty} \one_{A_n} d\mu \le \liminf_{n\to\infty}\int \one_{A_n}d\mu \\
&=\liminf_{n\to\infty} \mu(A_n)\le \liminf_{n\to\infty}\mu^*(B_n)+\eps=\lim_{n\to\infty}\mu^*(B_n)+\eps.
\end{align*}

The other preparation is a separation property for points and compact sets  in Hausdorff spaces: For $x\notin K$ compact, there
is a {\it closed} neighbourhood $U$ of $x$ which is disjoint from $K$ (so that $x$ and $K$ are separated by the disjoint open sets
$\Int(U)$ and $U^c$). Indeed, for each $y\in K$, there are disjoint open neighbourhoods $V_y$ of $y$ and $U_y$ of $x$. Compactness
yields a finite set $E\subseteq K$ with $K\subseteq\bigcup_{y\in E} V_y=V$ and then $U=\bigcap_{y\in E} U_y$ is a neighbourhood of $x$ such that
$U\cap V=\emptyset$ and hence $\ol U \subseteq \ol{V^c}=V^c\subseteq K^c$.

\begin{theo}[tightness of Borel measures]\label{thm:Tightness}
Every finite measure $\mu$  on the Borel \salg of a Suslin space $X$ is tight.
\end{theo}

\begin{proof}
Every non-empty Borel set $A\subseteq X$ is itself Suslin by theorem \ref{thm:Suslinstability} so that theorem \ref{thm:Suslin representation}
yields a continuous surjective $f:\N^\N \to A$. For $B_{\alpha|n}=\{\beta\in\N^\N:\beta|n=\alpha|n\}$ we get a Suslin
representation $F_{\alpha|n}=f(B_{\alpha|n})$ of $A$.
 For the pointwise partial order on $\N^\N$ and
 $G_{\alpha|n}=\bigcup\{F_{\beta|n}: \beta\le \alpha\}$, we then have
\begin{align*}
\mu(A) &= \mu\left(\bigcup_{\alpha_1\in\N} F_{(\alpha_1)} \right) \le \mu^*\left(\bigcup_{\alpha_1\in\N} G_{(\alpha_1)} \right)
\\ &\le \lim_{\alpha_1\to\infty} \mu^*(G_{(\alpha_1)}).
\end{align*}
For fixed $\eps>0$, there is thus $\alpha_1\in\N$ with $\mu(A)-\eps < \mu^*(G_{\alpha_1})$.
Since $G_{\alpha_1}$ is the increasing union of $G_{(\alpha_1,\alpha_2)}$, the same argument yields $\alpha_2\in\N$
with $\mu(A)-\eps < \mu^*(G_{(\alpha_1,\alpha_2)})$. Recursively, we obtain a sequence $\alpha\in\N^\N$ with
$\mu(A)-\eps <\mu^*(G_{\alpha|n})$ for all $n\in\N$. The closures $\ol G_{\alpha|n}$ form a decreasing sequence of
Borel sets and the continuity from
above for the {\it finite} measure $\mu$ yields,
for $M=\bigcap_{n\in\N} \ol G_{\alpha|n}$, 
\[
\mu(A)-\eps\le \lim_{n\to\infty}\mu\left( \ol G_{\alpha|n}\right)= \mu(M).
\]
It remains to see that $M$ is contained in a compact subset of $A$. This will follow from $M\subseteq f(L)$ with $L=\{\beta\in\N^\N: \beta\le\alpha\}$ because
$f:\N^\N\to A$ is continuous and $L$ is compact by Tychonov's theorem  (so that $K=f(L)$ is again compact, we have in fact
$M=K$ but this is not needed). For $x\in M$, the separation property from above yields that it is enough to show that every closed neighbourhood
$U$ of $x$ meets $K$. The definition of $G_{\alpha|n}$ and $x\in \ol G_{\alpha|n}$ yield  $\beta^n\in L$ with $f(\beta^n)\in U$ and for a convergent subsequence
$\beta^{n(\ell)}\to \beta \in L$ we obtain 
\[
f(\beta)=\lim\limits_{\ell\to\infty}f(\beta^{n(\ell)}) \in \ol U\cap f(L). \qedhere
\]
\end{proof}

Applying the theorem to complements we also get that finite Borel measures on Suslin spaces are {\it outer openly approximable}, i.e., for every
Borel set $A$, we have
$\mu(A)=\inf\{\mu(U): A\subseteq U \text{ open}\}$.

\section*{Acknowledgement}
It was a great pleasure for me to accompany Lara's journey into the fantastic world of formal proof assistance with Lean. We hope that her work will be soon available in the Lean mathematical library.
 I am also indebted to Julia Huschens for her careful reading as a {\it human proof checker}.

\bibliographystyle{amsalpha}
\bibliography{Effros}

\end{document}